\documentclass[11pt]{article}

\usepackage{amssymb,amstext,amsthm,latexsym}
\usepackage
{amsmath}
\usepackage{amsfonts}

\usepackage{mathrsfs}
\usepackage[bbgreekl]{mathbbol}
\usepackage{url}
\usepackage{doi}
\usepackage{mparhack}

\theoremstyle{definition}

\theoremstyle{remark}

\numberwithin{equation}{section}

\input{63arx.sty}

\begin{document}

\title{On Russell typicality in Set Theory\thanks
{Partial support of RFBR grant 20-01-00670 acknowledged.}}

\author{Vladimir~Kanovei\thanks 
{IITP RAS, Moscow, Russia, {\tt kanovei@iitp.ru}.}
\and
Vassily~Lyubetsky\thanks 
{IITP RAS, Moscow, Russia, {\tt lyubetsk@iitp.ru}}
}


\date{}

\maketitle

\begin{abstract}
By Tzouvaras, a set is nontypical in the Russell sense, 
if it belongs to a countable ordinal definable
set. 
The class $\HNT$ of all hereditarily nontypical sets 
satisfies all axioms of $\ZF$ and the double inclusion
$\HOD\sq\HNT\sq\rV$ holds. 
Several questions about the nature of such sets, 
recently proposed by Tzouvaras, 
are solved in this paper. 
In particular, a model of $\ZFC$ is presented in which 
$\HOD\sneq\HNT\sneq\rV$, and another model of $\ZFC$ 
in which $\HNT$ does not satisfy the axiom of choice. 
\end{abstract}

\maketitle

\punk{Introduction}

One of the fundamental directions in modern set theory
is the study of important classes of sets in the 
set theoretic universe $\rV$, which themselves
satisfy the axioms of set theory.
G\"odel's class $\rL$ of all \rit{constructive} sets 
traditionally belongs to such classes, as well as the class
$\HOD$ of all \rit{hereditarily ordinal definable} sets,
see \cite{jechmill} or \cite[\S\,7A]{kl1}. 
Both $\rL$ and $\HOD$ are transitive classes of sets 
in which all the axioms of the $\ZFC$ set theory, 
with the axiom of choice $\AC$, are fulfilled
(even if the universe $\rV$ itself only satisfies $\ZF$   
without the axiom of choice).
These classes satisfy $\rL\sq\HOD\sq\rV$,
and as it was established in early works on
modern axiomatic set theory,
the class $\HOD$ can be strictly between the classes
$\rL\sq\rV$ in suitable
generic extensions of $\rL$. 

Recent studies have shown considerable
interest in other classes of sets based on
the key concept of ordinal definability, which also
satisfy set-theoretic axioms.
In particular, the classes of 
\rit{nontypical} and
\rit{hereditarily nontypical}
sets are considered, 
whose name Tzouvaras \cite{tz20,tz21}
connects with philosophical and mathematical studies of 
Bertrand Russell and the works of van Lambalgen \cite{vL}
et al.\ on the axiomatization of the concept of randomness.

\bdf
\lam{nt}
The set $x$ is  
\rit{nontypical}, for short $x\in\NT$,
if it belongs to a countable 
\OD\ (ordinal definable) set.
The set $x$ is  
\rit{hereditarily nontypical}, for short $x\in\HNT$,
if it itself, all its elements, elements of elements,
and so on, are all nontypical, in other words
the \rit{transitive closure} $\tc(x)$ satisfies
$\tc(x)\sq\NT$.
\edf

The classes $\NT$ and $\HNT$ in this
definition correspond to $\NT_\ali$ and $\HNT_\ali$
in the basic definition system of \cite{tz21}.
Similarly defined narrower classes
$\NT_\alo$
(elements of finite ordinally definable sets)
and $\HNT_\alo$ in \cite{tz21} are identical to 
\rit{algebraically definable} and
\rit{hereditarily algebraically definable}
sets that have been investigated in recent papers
\cite{FGH,GH,HL}
and are not considered in this article. 

The class $\NT$ is not necessarily transitive, but the 
smaller class $\HNT\sq\NT$ is transitive and, as shown
in \cite{tz21}, satisfies all axioms of 
$\ZF$ (the axiom of choice \AC\ not included), 
and also satisfies
the relation $\HOD\sq\HNT\sq\rV$.
As for the axiom of choice $\AC$, if
$\rV=\rL$ (the constructibility axiom), then  
$\HNT=\HOD=\rL$ obviously holds, so in
this case $\AC$ holds in $\HNT$.

\bqe
[essentially Tzouvaras {\cite[\S\,2]{tz21}}]
\lam{q1}
Is it compatible with $\ZFC$ that
the axiom of choice $\AC$ {\ubf does not hold} in $\HNT$?
\eqe

The next problem of Tzouvaras {\cite[2.15]{tz21}}
aims to clarify the possibility of the precise 
\rit{equalities} in the relation $\HOD\sq\HNT\sq\rV$.

\bqe
\lam{q2}
Are the  next
sentences compatible with
\zfc\,?\vhm
\ben
\Renu
\itlb{q23}
$\HOD=\HNT\sneq\rV\,;$

\itlb{q22}
$\HOD\sneq\HNT=\rV\,;$

\itlb{q21}
$\HOD\sneq\HNT\sneq\rV\,.$
\een
\eqe

We answer all these questions in the positive. 
This is the main result of this article.
It is contained in the theorems
\ref{t3},\ref{t2},\ref{t4},\ref{t1} 
below.
The answer will be given through the construction 
of four corresponding
models of $\ZFC$ by the method of generic extensions
of the constructible universe $\rL$.

We begin with a model for Problem \ref{q2}\ref{q23} 
in Section \ref{x3}. 
It occurs that it is true in the extension $\rL[a]$ 
of $\rL$ by a single Cohen generic real that 
$\rL=\HOD=\HNT\sneq\rL[a]$ (Theorem~\ref{t3} below). 
This is based on our earlier result \cite{kl31} that 
the Cohen real $a$ does not belong to a countable \OD\ 
set in $\rL[a]$. 

As for Problem \ref{q2}\ref{q22}, we make use 
(Section \ref{x2}) of a forcing notion $\bP$ 
introduced in \cite{kl22} in order to define a 
generic real $a\in\dn$ whose \dd{\Eo}equivalence class 
$\eko a$ is a lightface $\ip12$ set with no \OD\ element. 

A positive answer to Problem \ref{q2}\ref{q21} is given 
in Sections \ref{x4}--\ref{x6} by means of the forcing 
notion $\bP\ti\text{(the Cohen forcing)}$. 
This includes the study of some aspects of the behavior 
of Borel functions in Sections \ref{bnf},\ref{nf}.

The article ends with a positive solution to 
Problem \ref{q1} in Section \ref{x1}. 
We make use of the finite-support product $\dP\lom$ 
of Jensen's forcing notion as in \cite{jenmin}.

\punk{Model I in which not all sets are nontypical}
\las{x3}

The following theorem solves the problem 
\ref{q2}\ref{q23} in the positive.
We make use of a well-known forcing notion.

\bte
\lam{t3}
If \/$a\in\dn$ is a Cohen generic real
over\/ $\rL$ then it is true in\/ $\rL[a]$ that\/
$\rL=\HOD=\HNT\sneq\rL[a]$.
\ete

Recall that Cohen generic extensions involve 
the forcing notion $\dC=2^{<\om}$
(all finite dyadic sequences).
Countable \OD\ sets in Cohen extensions
are investigated in our
papers \cite{kl39,kl31,kl40}.
In particular, we'll use
the following result here. 

\ble
[Thm 1.1 in {\cite{kl31}}]
\lam{l3}
Let \/$a\in\dn$ be Cohen-generic 
over the set universe\/ $\rV$.
Then it holds in\/ $\rV[a]$ that if\/
$Z\sq\dn$ is a countable \OD\ set then\/
$Z\in\rV.$\qed
\ele

This result admits the following extension 
for the case  $\rV=\rL$:

\bcor
\lam{c3}
Let \/$a\in\dn$ be Cohen-generic 
over the constructive universe\/ $\rL$.
Then it holds in\/ $\rL[a]$ that
if\/ $X\in\rL$ and\/ $A\sq2^X$ is 
countable\/ \OD\ then \/$A\sq\rL$.
\ecor
\bpf
As $\dC$ is countable,
there is a set of $Y\sq X\yt Y\in\rL$, 
countable in $\rL$ and such that if $a\ne b$
belong to $2^X$ then $a(x)\ne b(x)$ for some
$x\in Y$.
Then $Y$ is countable and \OD\ in $\rL[a]$,
so the \rit{projection}
$B=\ens{a\res Y}{a\in A}$ of the set $A$ will also
be countable and \OD\ in $\rL[a]$.
We have $B\sq\rL$ by the lemma.
(The set $Y$ here can be identified with $\om$.)
Hence, each $b\in B$ is \OD\ in $\rL[a]$.
However, if $a\in A$ and $b=a\res Y$,
then by the choice of $Y$ it holds in
$\rL[a]$ that 
$a$ is the only element in $A$ satisfying 
$a\res Y=b$.
Hence, $a\in\OD$.
\epf

\bpf[Theorem \ref{t3}]
The fact that $\rL=\HOD$ in $\rL[a]$ is a standard
consequence of the homogeneity of the 
Cohen forcing $\dC$.
Further, it is clear that $\HOD\sq\HNT$.
Let's prove the inverse relation 
${x\in\HNT}\imp{x\in\rL}$ 
in $\rL[a]$ by induction on
the set-theoretic rank $\rk x$ of 
sets $x\in\rL[a]$.
Since each set consists only of sets
of strictly lower rank, it is sufficient to check that
if a set $H\in\rL[a]$ satisfies
$H\sq\rL$ and $H\in\HNT$ in $\rL[a]$ then  
$H\in\rL$.
Here we can assume that in fact $H\sq\Ord$,
since $\rL$ allows an \OD\ wellordering.
Thus, let $H\sq\la\in\Ord$.
Additionally, since $H\in\HNT$,
we have, in $\rL[a]$, a countable \OD\ set
$A\sq\pws\la$ containing $H$.
However, $A\in\rL$ by Corollary \ref{c3}.
This implies $H\in\rL$. 
\epf

\punk{Perfect trees and \sit{s}}
\las{sst}

Our results will involve forcing notions 
that consist of perfect trees and \sit{s}.
Here we introduce the relevant
terminology from our earlier
works \cite{kl22,kl34,kl30}.

By $\bse$ we denote the set of all \rit{tuples}
(finite sequences) of terms $0,1$,
including the empty tuple $\La$.
The length of a tuple $s$ is denoted by $\lh s$,
and
$2^n=\ens{s\in\bse}{\lh s=n}$ 
(all tuples of length $n$).
A tree $\pu\ne T\sq\bse$ is \rit{perfect},
symbolically $T\in\pet$,
if it has no endpoints and isolated branches.
In this case, the set
$$
[T]=\ens{a\in\dn}{\kaz n\,(a\res n\in T)}
$$
of all \rit{branches} of $T$  
is a perfect set in $\dn.$
Note that $[S]\cap[T]=\pu$ iff $S\cap T$ is finite.\vhm

\bit
\item
If $u\in T\in\pet$, then
a \rit{portion}
(or a \rit{pruned tree}) 
$T\ret u\in\pet$ is defined by 
$T\ret u=\ens{s\in T}{u\su s\lor s\sq u}$.\vhm

\item
A tree $S\sq T$ is \rit{clopen} in $T$ iff it is 
equal to the union of a finite number of portions of $T$. 
This is equivalent to $[S]$ being 
clopen in $[T]$.\vhm
\eit

A tree\/ $T\sq\bse$ is a {\em Silver tree\/},
symbolically\/ $T\in\pes$, if 
there is an infinite
sequence of tuples $u_k=\uu kT\in\bse,$
such that $T$ consists of all tuples of the form
$$
s=u_0\we i_0\we u_1\we i_1\we u_2\we i_2
\we\dots\we u_{n}\we i_n
$$
and their sub-tuples, where $n<\om$ and $i_k=0,1$.
In this case the \rit{stem} $\stem T=\uu0T$ is equal 
to the largest tuple $s\in T$ with
$T=T\ret s$, and  
$[T]$ consists of all infinite sequences
$a=u_0\we i_0\we u_1\we i_1\we u_2\we i_2 \we\dots\in\dn,$
where $i_k=0,1$, $\kaz k$.
Put  
$$
\oin T{n}=\lh{u_0}+1+\lh{u_1}+1+\dots+\lh{u_{n-1}}
+1+\lh{u_n}\,.
$$
In particular, $\oin T{0}=\lh{u_0}$.
Thus 
$\oi T=\ens{\oin Tn}{n<\om}\sq\om$ is the set of all
\rit{splitting levels} of the Silver tree $T$. 

{\ubf Action.} 
Let $\sg\in\bse.$  
If $v\in\bse$ is another tuple of length
$\lh v\ge \lh\sg$, 
then the tuple $v'=\sg\aq v$ of the same length
$\lh {v'}=\lh v$ is defined by 
$v'(i)=v(i)+_2 \sg(i)$
(addition modulo $2$) for all $i<\lh\sg$,
but $v'(i)=v(i)$ whenever $\lh\sg\le i<\lh v$.
If $\lh v<\lh\sg$,
then we just define $\sg\aq v=(\sg\res\lh v)\aq v$.

If $a\in\dn,$ then similarly $a'=\sg\aq a\in\dn,$
$a'(i)=a(i)+_2\sg(i)$ for $i<\lh\sg$,
but $a'(i)=a(i)$ for $i\ge \lh\sg$.
If $T\sq\bse\yt X\sq\dn$, then the sets 
$$
\sg\aq T=\ens{\sg\aq v}{v\in T}
\qand
\sg\aq X=\ens{\sg\aq a}{a\in X}
$$
are \rit{shifts} of the tree $T$ and the set $X$  
accordingly.

\ble
[\cite{kl30}, 3.4]
\lam{LL}
If\/ $n< \om$ and\/ $u,v\in T \cap 2^n,$ then\/
$T \ret u=v \aq u \aq (T\ret v)$.

If\/ $t\in T\in\pes$ and\/ $\sg\in\bse,$ then\/
$\sg\aq T\in\pes$ and\/ $T\ret s\in\pes$.\qed
\ele

\vyk{
{\ubf Splittings.} 
\rit{The simple splitting} of a tree $T\in\pes$ consists of
subtrees
$$
\raw T0=T\ret{\roo T\we 0}
\qand 
\raw T1=T\ret{\roo T\we 1}\,,
$$
so $\raw Ti\in\pes$.
Splittings can be iterated: if
$s\in 2^n$, then
$$
\raw Ts=\raw{\raw{\raw
{\raw T{s(0)}}{s(1)}}{s(2)}\dots}{s(n-1)}
\in\pes
\,.
$$
Separately, put $\raw T\La=T$, for the empty tuple $\La$.

\ble
[\cite{kl30}, Lemma 4.2]
\lam{sadd}
Let\/ $T\in\pes\yt n<\om\yt h=\oin Tn.$
Then
\ben
\renu
\itlb{sadd0}%
if\/ $s\in 2^n$ then\/ $\raw Ts\in\pes$,
$\lh{(\roo{\raw Ts})}=h$,
and the tuple\/
$u[s]=\roo{\raw Ts}\in 2^h\cap T$ 
satisfies\/ $\raw Ts=\req T{u[s]}\;;$

\itlb{sad4+}%
hence,
$\ens{\raw Ts}{s\in 2^n}=\ens{\req Tu}{u\in 2^h\cap T}\;;$

\itlb{sad1}%
$T=\bigcup_{s\in 2^n}\raw Ts$,
and if\/ $s\ne t\in 2^n$ then\/
$[\raw Ts] \cap[\raw Tt]=\pu.$\qed
\een
\ele
}

\bdf
[\ubf refinements] 
\lam{dr}
Assume that $T,S\in\pes$, $S\sq T$, $n<\om$. 
We define $S\nq n T$
(the tree $S$ \rit{\dd nrefines\/ $T$}) 
if $S\sq T$ and $\oin Tk=\oin Sk$ for all $k<n$. 
This is equivalent to ($S\sq T$ and) 
$\uu kS=\uu kT$ for all $k<n$, of course.
\edf

Then $S\nq 0 T$ is equivalent to $S\sq T$,
and $S\nq{n+1} T$ implies $S\nq n T$ (and $S\sq T$),
but if $n\ge 1$ then $S\nq n T$ is equivalent to
$\oin T{n-1}=\oin S{n-1}$.

\ble
\lam{tadd}
Assume that\/ 
$T,U\in\pes\yt n<\om\yt h>\oin T{n-1} \yt s_0\in 2^h\cap T$,
and\/ $U\sq\req T{s_0}$. 
Then there is a unique tree\/ $S\in\pes$ such that\/
$S\nq n T$ and\/ $\req{S}{s_0}=U.$ 

If in addition\/ $U$ is clopen in\/ $T$ then\/ $S$ 
is clopen in\/ $T$ as well.
\ele
\bpf[sketch]
Define a tree $S$ so that 
$S\cap 2^h=T\cap 2^h$, and if
$t\in T\cap2^h$ then, by Lemma \ref{LL},
$\req {S}t=(t\aq s_0)\aq U$;
then $\req{S}{s_0}=U$. 
To check that $S\in\pes$, we can easily compute the 
tuples $\uu kS$. 
Namely, as $U\sq\req T{s_0}$, we have 
$s_0\sq \uu0U=\roo U$, hence 
$\ell=\lh{(\uu 0U)}\ge h>m=\oin T{n-1}$. 
%
Then $\uu kS=\uu kT$ for all $k<n$, 
$\uu nS=\uu 0U\res[m,\ell)$ 
(thus $\uu nS\in 2^{\ell-m}$), 
and $\uu kS=\uu kU$ for all $k>n$.
%
\epf

\ble
[\cite{kl30}, Lemma 4.4]
\lam{fus}
Let\/
$\ldots \nq 4 T_3\nq 3 T_2\nq 2 T_1\nq 1 T_0$ be 
a sequence of trees in\/ $\pes$.
Then\/ 
$T=\bigcap_nt_n\in\pes$.\qed
\ele
\bpf[sketch]
By definition we have $\uu k{T_n}=\uu k{T_{n+1}}$ 
for all $k\le n$. 
Then one easily computes that $\uu nT=\uu n{T_n}$ for 
all $n$.
\epf

\punk{Model II
in which there are more nontypical sets than \HOD\ sets}
\las{x2}

The following theorem solves the problem 
\ref{q2}\ref{q22} positively.

\bte
\lam{t2}
There is a generic extension of the constructible
universe\/ $\rL$, in which it is true that\/
$\HOD\sneq\HNT=\rV$.
\ete

Recall that the equivalence relation $\Eo$
is defined on $\dn$ so that $a\Eo b$ iff the set
$a\sd b=\ens{k}{a(k)\ne b(k)}$ is finite. 

To prove Theorem \ref{t2}, we will use an \OD\  
\dd{\Eo}equivalence class 
$$
\eko a=\ens{b\in\dn}{a\Eo b}=\ens{\sg\aq a}{\sg\in\bse}
$$ 
of a non-\OD\ generic real $a\in\dn$,
introduced in \cite{kl22} and also applied in
\cite{kl25,kl30,kl34}. 
This is done by a forcing notion 
$\bP$ having the following key properties, 
see \cite{kl22}.\vhm
\ben
\senu
\itlb{r1}%
$\bP\in\rL$ consists of Silver trees:
$\bP\sq\pes$.\vhm

\itlb{r2}%
If $u\in T\in\bP$ and $\sg\in\bse$ then
$T\ret u\in\bP$
and $\sg\aq T\in\bP$ --- this is the 
property of \rit{invariance}
\poo\ shifts and portions.\vhm

\itlb{r3}%
$\bP$ satisfies the countable antichain condition CCC  
in $\rL$.\vhm

\itlb{r5}%
The forcing $\bP$ ajoins a generic real $a\in\dn$ 
to $\rL$, whose \dd\Eo class
$\eko a=\ens{b\in\dn}{b\Eo a}$
is a (countable) \OD, and even $\ip12$
(lightface) set in
$\rL[a]$.\vhm

\itlb{r6}%
If a real $a\in\dn$ is $\bP$-generic
over $\rL,$ then $a$ is not \OD\ in
the generic extension $\rL[a]$.
(This property is an elementary consequence
of the invariance property as in \ref{r2},
see Lemma 7.5 in \cite{kl22}.)
\een

\bpf[Theorem \ref{t2}]
Let a real $a\in\dn$ be
$\bP$-generic over $\rL$.
According to \ref{r5} the real $a$
itself belongs to $\HNT$ in $\rL[a]$, hence 
the equality $\HNT=\rV$ holds 
in $\rL[a]$.
On the other hand, $a\nin\OD$ in $\rl[a]$
by \ref{r6}, thus $\HOD\sneq\HNT$
in $\rL[a]$, as required.
(A more thorough analysis based on \ref{r2}
shows that $\HOD=\rL$ in $\rL[a]$.)
\epf

\punk{Model III: nontypical sets in general position}
\las{x4}

The following theorem positively solves Problem 
\ref{q2}\ref{q21}, providing a model in which
hereditarily nontypical sets are strictly between 
$\HOD$ and $\rV$.

\bte
\lam{t4}
There is a generic extension of\/  
$\rl$, in which\/ 
$\HOD\sneq\HNT\sneq\rV.$
\ete

The proof of this theorem 
(to be completed in Section \ref{x6})
is based on a combination of ideas from
the proof of theorems \ref{t3} and \ref{t2}.
In fact, the forcing notion involved 
will be equal to the product $\bP\ti\dC$.
However, we will have to consider in more detail the
inductive construction of the set $\bP$,
as well as some questions related to continuous
and Borel functions and the construction of \sit{s}.

In the remainder, if $v\in\nse$
(a tuple of natural numbers), then we define 
$\bi v=\ens{x\in\bn}{v\su x}$, 
the {\em Baire interval\/} or {\em portion\/} 
in the Baire space $\bn.$

\punk{Reduction of  Borel functions to continuous ones}
\las{bnf}

A classical theorem claims that in Polish
spaces every Borel function is continuous on
a suitable dense $\Gd$ set
(Theorem 8.38 in Kechris \cite{Kdst}). 
It is also known that a Borel map defined on $\dn$ 
is continuous on a suitable \sit. 
The next lemma combines these two results. 
Our interest in functions defined on $\dn\ti\bn$ is 
motivated by further applications to reals in 
generic extensions of the form $\rL[a,x]$, where 
$a\in\dn$ is \dd\bP generic real for some $\bP\sq\pes$ 
while $x\in\bn$ is just Cohen generic.

\ble
\lam{su}
Let\/  $T\in\pes$ and 
$f:\dn\ti \bn\to\dn$ be a Borel map.
There is a \sit\/ $S\sq T$ and
a dense\/ $\Gd$ set\/ $D\sq\bn$
such that\/ $f$ is continuous on\/ $[S]\ti D.$
\ele
\bpf
By the abovementioned classical theorem,
$f$ is already continuous on some dense
$\Gd$ set $Z\sq[T]\ti\bn.$
It remains to define a \sit\/ $S\sq T$ and
a dense $\Gd$ set
$D\sq\bn$ such that $[S]\ti D\sq Z.$
This will be our goal.
 
We have $Z=\bigcap_nZ_n$,
where each $Z_n\sq [T]\ti\bn$ is open dense.

We will define
$S=\bigcup_k S_k$, where \sit s $S_k$
satisfy  
$\ldots \nq 4 S_3\nq 3 S_2\nq 2 s_1\nq 1 S_0=T$
as in the lemma \ref{fus}.
Tuples $w_k\in\nse\yt k<\om$ will also be defined.

We fix a recursive enumeration
$\om\ti\nse=\ens{\ang{N_k,v_k}}{k<\om}$. 

At step $0$ we already have $S_0=T$.
Assume that the tree $S_k$ has already been defined.
We claim that there exist:\vhm
\ben
\Aenu
\itlb{step4}%
a tuple $w_k\in\nse$ and a \sit\ $S_{k+1}\nq{k+1}S_k$,
{\ubf clopen} in $S_k$ (see Section \ref{sst}),
such that $v_k\sq w_k$ and
$[S_{k+1}]\ti \bi{w_k}\sq Z_{N_k}$.\vhm
\een
Now let $N=N_k\yt v=v_k$. 
Put $h=\oin{S_k}{k+1}$.
Consider any tuple $t\in 2^{h}\cap S_k.$
Since $Z_N$ is open dense,
there exist a tuple 
$v_1\in\nse$ and a \sit\ $A\sq\req{S_k}t$,
clopen in $S_k$ (for example, a portion in $S_k$)
such that $v\sq v_1$ and
$[A]\ti\bi{v_1}\sq Z_N$.
According to Lemma \ref{tadd}, there exists
a \sit\ $U_1\nq{k+1}S_k$,
clopen in $S_k$ along with $A$,
such that $\req{U_1}t=A$, so
$[\req{U_1}t]\ti\bi{v_1}\sq Z_N$ 
by construction. 

Now take another tuple $t'\in 2^{h}\cap S_k,$
and similarly find 
$v_2\in\nse$ and a \sit\ $A\sq\req{U_1}{t'}$,
clopen in $U_1$,
such that $v_1\sq v_2$ and
$[A]\ti\bi{v_2}\sq Z_N$.
Once again there is a \sit\ $U_2\nq{k+1}U_1$,
clopen in $S_k$ and such that 
$[\req{U_2}{t'}]\ti\bi{v_2}\sq Z_N$.

We iterate this construction over all tuples 
$t\in 2^{h}\cap S_k,$
\dd{\nq{k+1}}shrinking trees and extending tuples
in $\nse.$
We get a \sit\
$U\nq{k+1}S_k$,
clopen in $S_k$, and tuple
$w\in\nse,$ that $v\sq w$ and
$[U]\ti\bi{w}\sq Z_N$.
Take $w_k=w\yt S_{k+1}=U$.
This completes the inductive step.

As a result we get a sequence 
$\ldots \nq 4 S_3\nq 3 S_2\nq 2 S_1 \nq 1 S_0=T$
of \sit{s} $S_k$,
and tuples $w_k\in\nse$ ($k<\om$),
which satisfy \ref{step4} for all $k$.

We put $S=\bigcap_k S_k$;
then $S\in\pes$
by Lemma \ref{fus}, and $S\sq T$.

If $n<\om$ then let $W_{n}=\ens{w_k}{N_k=n}$.
Then $D_{n}=\bigcup_{w\in W_{n}}\bi w$
is an open dense set in $\bn.$
Indeed, let $v\in\nse.$ 
Consider $k$
such that that $v_k=v$ and $N_k=n$.
By construction, we have $v\sq w_k\in W_{n}$,
as required.
We conclude that the set
$D=\bigcap_{n}D_{n}$ is dense $\Gd$.

To check  $[S]\ti D\sq Z$,
let $n<\om$; we show that $[S]\ti D\sq Z_n$.
Let $a\in[S]$ and $x\in D$, in particular
$x\in D_n$, so $x\in\bi{w_k}$ for some
$k$ with $N_k=n$.
However, $[S_{k+1}]\ti\bi{w_k}\sq Z_n$ by \ref{step4}, 
and at the same time obviously $a\in[S_{k+1}]$.
We conclude that in fact
$\ang{a,x}\in Z_n$, as required.
\epF{Lemma \ref{su}}

\punk{Normalization of Borel maps}
\las{nf}

\bdf
\lam{nor}
A map $f:\dn\ti\bn\to\dn$ is 
\rit{normalized on\/ $T\in\pes$ for\/ $\bU\sq\pes$}
if there exists a dense\/ $\Gd$ set\/
$X\sq\bn$ such that $f$ is continuous on 
$[T]\ti X$ and$:$\vom
\bit
\item[$-$]
either $(1)$
there are tuples 
$v\in\nse\yt \sg\in\bse$ 
such that\/ $f(a,x)=\sg\aq a$ for all 
$a\in [T]$ and $x\in\bi v\cap X$, where
$\bi v=\ens{x\in\bn}{v\su x}\,;$\vom

\item[$-$]
or $(2)$
$f(a,x)\nin\bigcup_{\sg\in\bse\land S\in \bU}\sg\aq[S]$ 
for all $a\in [T]$ and $x\in X.$\qed
\eit
\eDf

\bte
\lam{tn}
Let\/ $\bU=\ans{T_0,T_1,T_2,\dots}\sq\pes$ and\/
$f:\dn\ti\bn\to\dn $ be a Borel map.
There is a set\/
$\bU'=\ans{S_0,S_1,S_2,\dots}\sq\pes$,   
such that\/ $S_n\sq T_n$ for all\/ $n$ and\/
$f$ is normalized on\/ $S_0$ for\/ $\bU'$.
\ete

\bpf
First of all, according to Lemma \ref{su}, there
is a \sit\ $T'\sq T_0$ and a dense ${\Gd}$ set
$W\sq\bn$ such that $f$ is continuous on $[T']\ti W$.
And since any dense
$\Gd$ set $X\sq\bn$ is homeomorphic to $\bn,$
we can \noo\ assume  
that $W=\bn$ and $T'=T_0$.
Thus, we simply suppose that $f$ is already
\rit{continuous on $[T_0]\ti\bn.$}

Assume that option (1) of the definition of \ref{nor}
does not take place, \ie\vom
\ben
\fenu
\itlb{*}%
if $X\sq\bn$ is dense $\Gd$, and
$v\in\nse\yt \sg\in\bse\yt S\in\pes\yt S\sq T_0$, 
then there are reals
$a\in [S]$ and $x\in\bi v\cap X$ such that 
$f(a,x)\ne\sg\aq a$.\vom
\een
We'll define \sit s $S_n\sq T_n$ and a dense
$\Gd$ set $X\sq\bn$ satisfying (2) of 
Definition \ref{nor}, that is, in our case,
the relation
$f(a,x)\nin\bigcup_{\sg\in\bse\land n<\om}\sg\aq[S_n]$ 
will be fulfilled 
for all $a\in [S_0]$ and $x\in X.$ 

The construction of the trees is organized in the form
$S_n=\bigcup_k S^n_k$, where \sit s $S^n_k$
satisfy  
$\dots \nq 4 S^n_3\nq 3 S^n_2\nq 2 S^n_1\nq 1 S^n_0=T_n$ 
as in Lemma \ref{fus} for each $n<\om$.
A series of tuples $w_k\in\nse$ ($k<\om$) will also be defined,
they will help us to construct a dense $\Gd$ set $X\sq\bn$ 
required.

We fix any enumeration
$\om\ti\bse\ti\nse
=\ens{\ang{N_k,\sg_k,v_k}}{k<\om}$. 

At step $0$  of the construction, we put 
$S^n_0=T_n$ for all $n$.

Assume that $k<\om$ and all \sit s
$S^n_k\yt n<\om$ are already defined.\vhm
\ben
\Aenu
\atc
\itlb{step}%
We put $S^n_{k+1}=S^n_k$ for all $n>0\yt n\ne N_k$.
\een
As for the trees $S^0_{k+1}$ and $S^{N_k}_{k+1}$,
we claim that there exist:\vhm
\ben
\Aenu
\atc
\atc
\itlb{step1}%
a tuple $w_k\in\nse$ and \sit s 
$S^0_{k+1}\nq{k+1}S^0_k\yt S^{N_k}_{k+1}\nq{k+1}S^{N_k}_k$  
such that $v_k\sq w_k$ and
$f(a,x)\nin\sg_k\aq[S^N_{k+1}]$ 
for all $a\in [S^0_{k+1}]$ and $x\in\bi {w_k}$.\vom
\een
For brevity, let $N=N_k\yt\sg=\sg_k\yt v=v_k$. 
Put $h=\oin{S^0_{k}}{k+1}$, 
$m=\oin{S^N_{k}}{k+1}$.\vom

\rit{Case 1\/}: $N>0$.
Take any pair of tuples $s\in2^{h}\cap S^0_{k}$, 
$t\in2^{m}\cap S^N_{k}$
and any reals $a_0\in[\req{S^0_k}s]$ and $x_0\in\bn.$
Consider any real $b_0\in[\req{S^N_k}t]$ not equal
to $\sg\aq f(a_0,x_0)$.
Let's say $b_0(\ell)=i\ne j= (\sg\aq f(a_0,x_0))(\ell)$,
where $i,j\le1\yt \ell<\om$.
By the continuity of $f$, there is a tuple 
$v_1\in\nse$ and \sit\ $A\sq\req{S^0_k}s$
such that $v\sq v_1\su x_0$, $a_0\in[A]$, and
$(\sg\aq f(a,x))(\ell)=j$ for all $x\in\bi{v_1}$
and $a\in [A]$.
It is also clear that
$B=\ens{\tau\in \req{S^N_k}t}
{\lh \tau\le\ell\lor \tau(\ell)=i}$ 
is a \sit\ containing $b_0$, 
and $b(\ell)=i$ for all $b\in[B]$.
According to Lemma \ref{tadd}, there are  
\sit{s} $U_1\nq{k+1}S^0_k$ and $V_1\nq{k+1}S^N_k$,
such that $\req{U_1}s=A$ and $\req{V_1}t=B$, hence 
by construction we have
$\sg\aq f(a,x)\nin [\req{V_1}t]$ for all
$a\in[\req{U_1}s]$ and $x\in\bi{v_1}$.

Now consider another pair of tuples 
$s\in2^{h}\cap S^0_{k}$, $t\in2^{m}\cap S^N_{k}$.
We similarly get Silver trees
$U_2\nq{k+1} U_1$ and $V_2\nq{k+1}V_1$, and a tuple
$v_2\in\nse,$ such that $v_1\sq v_2$ and
$\sg\aq f(a,x)\nin [\raw{V_2}{t'}]$ for all
$a\in[\req{U_2}{s'}]$ and $x\in\bi{v_2}$.
In this case, we have $\req{V_2}{t}\sq\req{V_1}{t}$
and $\req{U_2}{s}\sq\req{U_1}{s}$, so that what has already
been achieved at the previous step is preserved.

We iterate through all pairs of 
$s\in2^{h}\cap S^0_{k}$, $t\in2^{m}\cap S^N_{k}$,
\dd{\nq{k+1}}shrinking trees and extending tuples
in $\nse$ at each step.
This results in a pair of Silver trees
$U\nq{k+1}S^0_k\yt V\nq{k+1}S^N_k$ and a tuple
$w\in\nse$ such that $v\sq w$ and
$\sg\hspace{0.1ex}\aq f(a,x)\nin [V]$ for all
reals $a\in[U]$ and $x\in\bi{w}$.
Now to fulfill \ref{step1}, take $w_k=w$,
$S^0_{k+1}=U,$  and $S^{N_k}_{k+1}=V.$ 
Recall that here $N_k=N>0$.\vom

\rit{Case 2\/}: $N=0$.
Here the construction somewhat changes, and
hypothesis \ref{*} will be used. 
We claim that there exist:
\ben
\Aenu
\atc
\atc
\atc
\itlb{step0}%
a tuple $w_k\in\nse$ and a Silver tree
$S^0_{k+1}\nq{k+1}S^0_k$
such that $v_k\sq w_k$ and
$f(a,x)\nin\sg_k\aq[S^0_{k+1}]$ 
for all $a\in [S^0_{k+1}]$ and $x\in\bi {w_k}$.
\een
As above, let $\sg=\sg_k\yt v=v_k$.
Take any pair of tuples
$s,t\in2^{h}\cap S^0_{k}$, where  
$h=\oin{S^0_k}{k+1}$ as above.
%
%
Thus 
$S^0_k\ret{t}=t\aq s\aq (S^0_k\ret{s})$,
by Lemma~\ref{LL}.
According to \ref{*}, there are reals $x_0\in\bi{v}$ and
$a_0\in [S^0_k\ret{s}]$ satisfying
$f(a_0,x_0)\ne \sg\aq s\aq t\aq a_0$, or 
equivalently, 
$\sg\aq f(a_0,x_0)\ne s\aq t\aq a_0$. 

Similarly to Case 1, we have
$(\sg\aq f(a_0,x_0))(\ell)=i\ne 
j= (s\aq t\aq a_0)(\ell)$ 
for some $\ell<\om$ and $i,j\le 1$.
By the continuity of $f$, there is a tuple 
$v_1\in\nse$ and a \sit\  
$A\sq S^0_k\ret s$, clopen in $S^0_k$, 
such that $v\sq v_1\su x_0$, $a_0\in[A]$, and
$(\sg\aq f(a,x))(\ell)=j$ but 
$(s\aq t\aq a)(\ell)=j$  
for all $x\in\bi{v_1}$ and $a\in [A]$.
Lemma \ref{tadd} gives us a 
\sit\ $U_1\nq{k+1}S^0_k$, clopen in $S^0_k$ as well,  
such that $\req{U_1}s=A$ ---
and then $\req{U_1}t=s\aq t\aq A$.
Therefore 
$\sg\aq f(a,x)\nin [\req{U_1}t]$ 
holds for all
$a\in[\req{U_1}s]$ and $x\in\bi{v_1}$ 
by construction. 

Having worked out all pairs of tuples 
$s,t\in2^{h}\cap S^0_{k}$,
we obtain a \sit\ $U\nq{k+1}S^0_k$ and a tuple
$w\in\nse,$ such that $v\sq w$ and
$\sg\aq f(a,x)\nin[U]$ for all
$a\in[U]$ and $x\in\bi{w}$.
Now to fulfill \ref{step0}, take $w_k=w$
and $S^0_{k+1}=U$.

To conclude, we have for each $n$ a sequence
$\ldots \nq 4 S^n_3\nq 3 S^n_2\nq 2 S^n_1\nq 1 S^n_0=T_n$
of \sit{s} $S^n_k$,
along with tuples $w_k\in\nse$ ($k<\om$),
and these sequences satisfy the requirements
\ref{step} and \ref{step1}
(equivalent to \ref{step0} in case $N_k=0$).

We put $S_n=\bigcap_k S^n_k$. 
Then $S_n\in\pes$
by Lemma \ref{fus}, and $S_n\sq T_n$.

If $n<\om$ then let
$W_{n\sg}=\ens{w_k}{N_k=n\land \sg_k=\sg}$.  
The set $X_{n\sg}=\bigcup_{w\in W_{n\sg}}\bi w$ 
is then open dense in $\bn.$
Indeed, if $v\in\bn$ then we take $k$
such that $v_k=v\yt N_k=n\yt\sg_k=\sg$;
then $v\sq w_k\in W_{n\sg}$ by construction.
Therefore,
$X=\bigcap_{n<\om\yi \sg\in\bse}X_{n\sg}$  
is a dense $\Gd$ set.
Now to check property (2)
of Definition \ref{nor}, consider 
any $n<\om\yt\sg\in\bse\yt a\in[S_0]\yt x\in X$;
we show that $f(a,x)\nin\sg\aq[S_n]$.

By construction, we have $x\in X_{n\sg}$, \ie\
$x\in\bi{w_k}$, where $k\in W_{n\sg}$, so that
$N_k=n\yt \sg_k=\sg$.
Now $f(a,x)\nin\sg\aq[S_n]$ directly
follows from \ref{step1} for this $k$, since
$S_0\sq S^0_{k+1}$ and $S_n\sq S^n_{k+1}$.
\epF{Theorem~\ref{tn}}

\punk{The  forcing notion for Model III}
\las{x5}

Using the standard encoding of Borel sets,
as \eg\ in \cite{sol}
or \cite[\S\,1D]{kl1}, we fix
a coding of Borel functions $f:\dn\ti\bn\to \dn.$
As usual, it includes a 
\rit{\dd{\ip11}set\snos
{The letters $\is{}{}$ and $\ip{}{}$ denote
effective (lightface) projective classes.}
of codes\/} 
$\bc\sq\bn$,
and for each code $r\in\bc$ a certain 
Borel function $F_r:\dn\ti\bn\to \dn$
coded by $r.$
We assume that each Borel function has some
code, and there is a 
$\is11$ relation $\gS(\cdot,\cdot,\cdot,\cdot)$ and
a $\ip11$ relation $\gP(\cdot,\cdot,\cdot,\cdot)$
such that for all $r\in\bc\yt x\in\bn,$ and $a,b\in\dn$
it holds 
$F_r(a,x)=b\eqv\gS(r,a,x,b)\eqv\gP(r,a,x,b)$.

If $\bU\sq\pes$, then $\clo\bU$ denotes
the set of all trees of the form $\sg\aq(T\ret s)$,
where $\sg\in\bse$ and $s\in T\in\bU$, \ie\ the closure
of $\bU$ \poo\ both shifts and portions.\vhm

{\ubf The following construction is maintained in $\rL$.}
We define a sequence of countable sets
$\bU_\al\sq\pes\yt\al<\omi$ satisfying the following 
conditions \ref{b1}--\ref{b5}.\vhm
\ben
\benu
\itlb{b1}%
Each $\bU_\al\sq\pes$ is countable,
$\bU_0$ consists of a single
tree $\bse.$\vhm
\een
We then define 
$\bP_\al=\clo{\bU_\al}$,
$\bpl\al=\bigcup_{\xi<\al}\bP_\xi$. 
These
sets are obviously closed with respect
to shifts and portions, that is  
$\clo{\bP_\al}={\bP_\al}$ and
$\clo{\bpl\al}={\bpl\al}$.\vhm
\ben
\benu
\atc
\itlb{b2}%
For every $T\in\bpl\al$ there is a tree  
$S\in \bU_\al\yt S\sq T$.\vhm
\een
Let $\zfcm$ be the subtheory of the theory \ZFC,
containing all axioms except the power set axiom, and
additionally containing an axiom asserting the existence
of the power set $\pws\om$.
This implies the existence of $\pws X$ for any
countable $X$, the existence of $\omi$ and $\dn$, as well
as the existence of continual sets like $\dn$ or $\pes$.

By $\gM_\al$ we denote the smallest model of 
$\zfcm$ of the form $\rL_\la$
containing the sequence
$\sis{\bU_\xi}{\xi<\al}$, in which $\al$ and all
sets $\bU_\xi\yt\xi<\al$, are countable.\vhm
\ben
\benu
\atc
\atc
\itlb{b3}%
If a set $D\in\gM_\al\yt D\sq \bpl\al$
is dense in $\bpl\al$, and $U\in\bU_\al$, then
$U\sqf D$,
meaning that there is a finite set $D'\sq D$ such that
$U\sq\bigcup D'$.\vhm
\een
Given that $\clo{\bpl\al}={\bpl\al}$,
this is automatically transferred to all
trees $U\in\bP_\al$ as well.
It follows that $D$ remains
predense in $\bpl\al\cup\bP_\al$.

To formulate the next property, we fix
an enumeration
$\pes\ti\bn=\ens{\ang{T_\xi,b_\xi}}{\xi<\omi}$ 
in $\rL$,
which 1) is definable in $\rL_{\omi}$, and
2) each value in $\pes\ti\bn$ is taken
uncountably many times.\vhm
\ben
\benu
\atc
\atc
\atc
\itlb{b4}%
If $T_\al\in\bpl\al$ then there is a tree $S\in\bU_\al$
satisfying $S\sq T$, on which $F_{b_\al}$
is normalized for $\bU_\al$ 
in the sense of Definition \ref{nor}.\vhm

\itlb{b5}%
The sequence $\sis{\bU_\al}{\al<\omi}$
is \dd\in definable in $\rL_\omi$.
\een

The construction goes on as follows. 
{\ubf Arguing in $\rL$},
suppose that $\al<\omi$,
the subsequence $\sis{\bU_\xi}{\xi<\al}$
has been defined, and the sets 
$\bP_\xi=\clo{\bU_\xi}$ (for $\xi<\al$), $\bpl\al$,
$\gM_\al$ are defined as above.

\ble
[in $\rL$]
\lam{susU}
Under these assumptions, there is a countable set\/
$\bU_\al\sq\pes$ satisfying\/
\ref{b2}, \ref{b3}, \ref{b4}.
\ele
\bpf
The existence of a countable set
$\bU_\al\sq\pes$ satisfying \ref{b2}, \ref{b3}
is known from our earlier papers, see  
\cite[\S\,4]{kl22}, \cite[\S\,9 and 10]{kl30},
\cite[\S\,10]{kl34}.
If now the tree $T_\al$ belongs to $\bpl\al$
(if not then we don't worry about \ref{b4}),
then we consider, according to \ref{b2}, a tree
$T\in\bU_\al$ satisfying $T\sq T_\al$.
Using Theorem \ref{tn}, we shrink each
tree $U\in\bU_\al$ to a tree $U'\in\pes\yt U'\sq U$,
so that the function $F_{b_\al}$ is normalized
on $T'$ for $\bU'=\ens{U'}{U\in\bU_\al}$.
Finally take $\bU'$ as 
$\bU_\al$ and $T'$ as $S$ to fulfill \ref{b4}.
\epF{Lemma}

To accomplish the construction, we take $\bU_\al$
to be the smallest, in the sense of the G\"odel 
wellordering of $\rL$,
of those sets that exist by Lemma \ref{susU}.
Since the whole construction is relativized to
$\rL_\omi$, the requirement \ref{b5} is also met.

We put $\bP_\al=\clo{\bU_\al}$ for all $\al<\omi$,
and $\bP=\bigcup_{\al<\omi} \bP_\al$.

The following result, in part related to CCC, 
is a fairly standard
consequence of \ref{b3}, see for example
\cite[6.5]{kl22}, \cite[12.4]{kl34}, or
\cite[Lemma 6]{jenmin}; we will skip the proof.

\ble
[in $\rL$]
\lam{cccU}
The forcing notion\/ $\bP$ belongs to\/ $\rL$, 
satisfies\/ $\bP=\clo{\bP}$ 
and satisfies CCC in\/ $\rL$.\qed
\ele

\ble
[in $\rL$]
\lam{norU}
Let\/ $T\in\bP$ and \/$f:\dn\ti\bn\to\dn $ be a 
Borel function.
There is an ordinal\/ $\al<\omi$ and a tree\/
$S\in\bU_\al$, $S\sq T$, on which\/ $f$
is normalized for\/ $\bU_\al$.
\ele

\bpf
By the choice of the enumeration 
of pairs in $\pes\ti\bn,$
there is an ordinal $\al<\omi$ such that 
$T\in\bpl\al$ and $T=T_\al$, $f=F_{r_\al}$.
It remains to refer to \ref{b4}.
\epf

\punk{Model III: finalization}
\las{x6}

We use the product $\bP \ti\dC$ of the 
forcing notion $\bP$ defined in $\rL$ 
and satisfying conditions \ref{b1}--\ref{b5}
as above, and the Cohen forcing,
here in the form of $\dC=\nse$.

\bte
\lam{t4.}
Let a pair of reals\/ $\ang{a_0,x_0}$ be\/
\dd{\bP\ti\dC}generic over \/ $\rL$.
Then it is true in\/ $\rL[a,x]$ that\/
$\HOD\sneq\HNT\sneq\rV,$ and more precisely$:$
\ben
\renu
\itlb{t4a}%
$a_0$ is not $\OD$ in $\rL[a_0,x_0]\;;$

\itlb{t4b}%
$a_0$ belongs to $\HNT$ in $\rL[a_0,x_0]\;;$

\itlb{t4c}%
$x_0$ does not belong to $\HNT$ in $\rL[a_0,x_0]\;.$
\een
\ete
\bpf
\ref{t4a}
By the forcing product theorem, $a_0$ is a 
\dd\bP generic real over $\rl[x]$.
However the forcing notion 
$\bP$ is invariant
\poo\ shifts by construction, that is if $T\in\bP$
and $\sg\in\bse$ then $\sg\aq T\in\bP$.
Now the result required is obtained
by an elementary argument,
see Lemma 7.5 in \cite{kl22}.

\ref{t4b}
It is sufficient to prove that the 
\dd\Eo equivalence class 
$\eko{a_0}$ of our generic
real $a_0$ is an \OD\ set in $\rL[a_0,x_0]$.
According to \ref{b5}, it suffices to establish
the equality
$$
\textstyle
\eko{a_0}=\bigcap_{\xi<\omi}\bigcup_{T\in\bP_\xi}[T]\,.
\eqno(*)
$$
Note that every set $\bP_\xi$ is pre-dense
in $\bP$; this follows from \ref{b3} and \ref{b4}, see,
for example, Lemma 6.3 in \cite{kl22}.
This immediately implies $a_0\in \bigcup_{T\in\bP_\xi}[T]$
for each $\xi$.
Yet all sets $\bP_\xi$ are invariant \poo\ 
shifts by construction.
Thus we have $\sq$ in (*).

To prove the inverse inclusion, assume
that a real $b\in\dn$ belongs to the right-hand side 
of (*) in $\rL[a_0,x_0]$.
It follows from Lemma \ref{cccU} 
(and the countability of $\DC$) 
that the forcing $\bP\ti\dC$ preserves cardinals.
We conclude that that $b=g(a_0,x_0)$ for some Borel
function $g=F_q:\dn\ti\bn\to\dn$ with a code $q\in\bc\cap\rL$.
{\ubf Assume to the contrary that} $b=g(a_0,x_0)\nin\eko{a_0}$.
Since $x_0\in\bn$ is a \dd\dC generic
real over $\rl[a_0]$ by the forcing product theorem,
this assumption is forced, so that there is a 
tuple $u\in\dC=\nse$ such that
$$
\textstyle
f(a_0,x)\in\bigcap_{\xi<\omi}\bigcup_{T\in\bP_\xi}[T]
\bez \eko{a_0}\,,
$$
whenever a real $x\in\bi u$ is
\dd\dC generic over $\rL[a_0]$.
(Recall that $\bi u=\ens{y\in\bn}{u\su y}$.)
Let $H$ be the canonical homomorphism of 
$\bn$ onto $\bi u$. 
We put $f(a,x)=g(a,H(x))$ for $a\in\dn\yt x\in\bn.$
Then $H$ preserves the \dd\dC genericity,
and hence
$$
\textstyle
f(a_0,x)\in\bigcap_{\xi<\omi}\bigcup_{T\in\bP_\xi}[T]
\bez \eko{a}\,,
\eqno(**)
$$
whenever $x\in\bn$ is 
\dd\dC generic over $\rL[a_0]$.
Note that $f$ also has a Borel code 
$r\in\bc$ in $\rL$, so that $f=F_{r}$.

It follows from Lemma \ref{norU} that there is an ordinal 
$\al<\omi$ and a tree
$S\in\bU_\al$, on which $f$ is normalized for $\bU_\al$, 
and which satisfies $a_0\in [S]$.
Normalization means that, in $\rL$,
there is a dense $\Gd$ set $X\sq\bn$
satisfying one of the two options of 
Definition~\ref{nor}.
Consider a real $z\in\bn\cap\rL$
(a \rit{$\Gd$-code} for $X$ in $\rl$) 
such that 
$X=X_z=\bigcap_k\bigcup_{z(2^k\cdot 3^j)=1}\bi{w_j}$,
where $\bse=\ens{w_j}{j<\om}$ is a 
fixed recursive enumeration of tuples.\vom

\rit{Case 1}:
there are tuples
$v\in\nse\yt \sg\in\bse,$ 
such that $f(a,x)=\sg\aq a$ for all 
points $a\in [S]$ and $x\in\bi v\cap X$.
In other words, it is true in $\rL$ that 
$$
\kaz a\in[S]\,\kaz x\in\bi v\cap X_z\:
(f(a,x)=\sg\aq a)\,.
$$
But this formula is absolute by Shoenfield, so it is also true
in $\rL[a_0,x_0]$.
Take $a=a_0$ (recall: $a_0\in[S]$)
and any real $x\in\bi v$,
\dd\dC generic over $\rL[a_0]$.
Then $x\in X_z$, because $X_z$ is a dense $\Gd$
with a code even from $\rL$.
Thus $f(a_0,x)=\sg\aq a_0\in\eko{a_0}$,
which contradicts (**).\vom

\rit{Case 2}:
$f(a,x)\nin\bigcup_{\sg\in\bse\land U\in \bU_\al}
\sg\aq[U]$ 
for all $a\in [S]$ and $x\in X.$
By the definition of $\bP_\al$, this implies 
$f(a,x)\nin\bigcup_{T\in \bP_\al}[T]$
for all $a\in [S]$ and $x\in X,$
and this again contradicts (**) for $a=a_0$.
\vom

The resulting contradiction in both cases refutes
the contrary assumption above and completes the 
proof \ref{t4b}.\vom

Finally, \ref{t4c} follows from Lemma \ref{l3}
for $\rV=\rL[a_0]$, since obviously 
$x_0\nin\rL[a_0]$.\vom

\epF{theorems \ref{t4.} and \ref{t4}}

\punk{Model IV: nontypical sets
sans the axiom of choice}
\las{x1}

The following theorem solves Problem 
\ref{q1} in the positive.

\bte
\lam{t1}
There is a generic extension of the constructible
universe\/ $\rL$ in which\/ $\AC$ holds 
but it is true that the class\/
$\HNT$ does not satisfy\/ $\AC$.
\ete

We will use a forcing notion $\dP\in\rL$ defined 
in \cite[\S\,7]{kl27} in order to obtain a model with 
a non-empty countable \OD\ 
set of pairwise generic reals, containing no 
\OD\ reals.
Modulo technical details,
this forcing coincides with the Jensen forcing
from \cite{jenmin} 
(also presented in \cite[28.A]{jechmill}).
The crucial step in \cite{kl27} was the proof
that those key properties of Jensen's forcing 
responsible for the uniqueness and definability 
of generic reals, previously established for $\dP$
and its finite products $\dP^n$, for example, in \cite{ena},
also hold for the countable product $\dP^\om.$
This forcing and its derivates were used in \cite{abr} 
and recently
in \cite{jml19,kl28,kl49,fm21} for various purposes.
This forcing $\dP$ has the following main  
properties \ref{p1}--\ref{p5}, see \cite{kl27}.

\ben
\cenu
\itlb{p1}%
$\dP\in\rL$, $\dP\sq\pet,$
$\dP$ contains the full tree $\bse.$

\itlb{p2}%
If $u\in T\in\dP$, then the \rit{portion}
$T\ret u$
also belongs to $\dP$.

\itlb{p3}%
$\dP$ satisfies CCC in $\rL$: each antichain 
$A\sq\dP$ is at most countable.

\itlb{p4}%
The set $\dP\lom,$ that is, 
the \rit{weak product}, or
\rit{product with finite support},
also satisfies CCC.
To be precise, here $\dP\lom$ consists of all
functions $\tau:\dom\tau\to\dP$, where
$\dom\tau\sq\om$ is finite.

\itlb{p5}%
Forcing $\dP\lom$ naturally adjoins a generic sequence 
of the form $\va=\sis{a_n}{n<\om}$ of 
$\dP$-generic reals $a_n\in\dn$ to $\rL.$
The corresponding set $\ww(\va)=\ens{a_n}{n<\om}$
is a (countable) \OD, and even $\ip12$
(without parameters) set in
the generic extension $\rL[\va]$.
\een

To prove Theorem \ref{t1}, we consider a 
$\dP\lom$-generic extension $\rL[\va]$
as in \ref{p5}, and 
the class $\HNT^{\rL[\va]}$ in this extension.
Our goal will be to prove that $\AC$ is false in
$\HNT^{\rL[\va]}$.
This will be a simple consequence of the last statement
of the next lemma.
In the remainder, if $W\sq\bn$
then $\qq W$ will denote Cohen forcing for adding
a generic $1-1$ function $f:\om\onto W$.
Thus, $\qq W$ consists of all $1-1$
functions $p:\dom p\to W$, where
$\dom p\sq\om$ is finite.

\ble
\lam{l1}
Let\/ $\va=\sis{a_n}{n<\om}$
be a\/ $\dP\lom$-generic sequence
over \/ $\rL$, and\/ $W=\ww(\va)$.
Then$\;:$\vhm
\ben
\renu
\itlb{l11}%
$\rL(W)\sq(\HOD)^{\rL[\va]}\;;$\vhm

\itlb{l12}%
$W$ is not a wellorderable set
in\/ $\rL(W)\;;$\vhm

\itlb{l13}%
$\va$ is a\/ $\qq W$-generic function
over\/ $\rL(W)\;;$\snos
{It is an important point here  that the
same function or sequence $\va\in W^\om$
can act as both a \dd{\dP\lom}generic object
over $\rL$
and as a \dd{\qq W}generic object over $\rL(W)$.
Moreover, the extensions $\rL[\va]$ and $\rL(W)[\va]$
coincide.
Such representations of a one-step generic
extension as a multi-step extension 
(here two-step) 
are well known, see, for example,
\cite{sol,gri}, \cite[\S\,7]{kl21}, \cite{kl32}.}\vhm

\itlb{l14}%
if\/ $\vb\in W^\om$ is a\/
$\qq W$-generic function over\/ $\rL(W)$
then\/ $\vb$ is a\/
$\dP\lom$-generic sequence over\/ $\rL$
in the sense of\/ \ref{p5} and\/ 
$\rL(W)[\vb]=\rL[\vb]\;;$\vhm

\itlb{l15}%
if\/ $\vb\in W^\om,$ the pair\/ $\ang{\va,\vb}$
is\/ \dd{(\qq W\ti\qq W)}generic over\/
$\rL(W),$ and\/ $Z\in\rL(W)[\va]\cap\rl(W)[\vb]$,
$Z\sq\rL(W)$, then \/$Z\in\rL(W)\;;$\vhm

\itlb{l16}%
if\/ $Z\in\rL[\va]$, $Z\sq W^\om$ is a countable\/
$\OD$ set in\/ $\rL[\va],$
then\/ $Z\sq\rL(W)\;.$
\een
\ele

\bpf[Theorem \ref{t1} from the lemma]
Here we show how Theorem \ref{t1} follows from 
the lemma,
and then we prove the lemma itself.
It suffices to prove that the set
$\ww=\ww(\va)=\ens{a_n}{n<\om}$, which belongs to
$\HNT^{\rL[\va]}$ according to \ref{p5}, is not  
wellorderable in $\HNT^{\rL[\va]}$.
Suppose to the contrary that such a wellordering exists.
Then there is also a bijection $f\in\HNT^{\rL[\va]}$,
$f:\om\na\ww$.
By definition, such a bijection belongs to a countable
\OD\ set $Z\in\rL[\va]$, $Z\sq W^\om,$ in $\rL[\va]$.
According to claim \ref{l16}, we have $Z\sq\rL(W)$,
so $f\in\rL(W)$, \ie\ $W$ is wellordered  
in $\rL(W)$, which gives a contradiction with 
claim \ref{l12} of the lemma.\vom

\epF{Theorem \ref{t1} from Lemma \ref{l1}}

\bpf[Lemma \ref{l1}] 
To prove \ref{l11}, note that $W$ is 
a countable $\OD$ set in $\rl[a]$ by \ref{p5},
therefore $\ww$ belongs to $\HNT$.

Further, \ref{l12} is a common property of permutation
models.

To prove \ref{l13}, assume towards the contrary 
that there is a set $D\in\rL(W)\yt D\sq\qq W$,
dense in $\qq W$, and such that no \usl{}\ $q\in D$ 
is extended by $\va$.
As an element of $\rL(W)$, the set $D$ is definable in
$\rL(W)$ in the form:
$$
D=\ens{q\in\qq W}{\vpi(q,W,a_0,\dots,a_n,x)}\,,
$$
where $x\in\rL$, $n<\om$, and $a_0,\dots,a_n$ are 
the initial terms of the sequence $\va$.
There is a \usl{}\ $\tau\in\dP\lom,$
which is compatible with $\va$ and 
\dd{\dP\lom}forces,
over $\rl$, our assumption. 
That is, if a $\dP\lom$-generic sequence
$\vb=\sis{b_n}{n<\om}$ extends $\tau$ then
the set $D(\vb)$ defined in $\rL(W(\vb))$ by
$$
D(\vb)=
\ens{q\in\qq{W(\vb)}}{\vpi(q,W(\vb),b_0,\dots,b_n,x)}\,,
$$
is dense in $\qq{W(\vb)}$, but no \usl{} $q\in D(\vb)$ 
is extended by $\vb$.

We can \noo\ assume that 
$\dom\tau=\ans{0,1,\dots,n}$. 

Now consider a \usl{} $p\in\qq W$ defined by 
$p(j)=a_j$ for all $j=0,1,\dots,n$. 
As $D$ is dense, there exists a \usl{} $q\in D$ 
extending $p$. 
Then $\dom q=\ans{0,1,\dots,n}\cup U$, where 
$U\sq\ans{n+1,n+2,\dots}$ is a finite set.
If $i\in U$ then by definiton  
$q(i)=a_{k_i}$, where $k_i \ge n+1$ and the map 
$i\mto k_i$ is injective. 

There is a bijection $\pi:\om\na\om$ satisfying 
$\pi(j)=j$ for all $j\le n$, $\pi(i)=k_i$ for  
all $i\in U$, and $\pi(\ell)=\ell$ generally for  
all but finite numbers $\ell<\om$, in particular, 
$\pi\in\rL$.
The sequence $\vb=\sis{b_n}{n<\om}$, defined 
by $b_i=a_{\pi(i)}$ for all $i<\om$, 
is $\dP\lom$-generic by the choice of $\pi$, 
and obviously $W(\vb)=W(\va)=W$. 
In addition, $b_j=a_j$ for all $j=0,1,\dots,n$, 
thus $\vb$ extends $\tau$.
We also have $D(\vb)=D(\va)=D$, and hence  
the abovedefined \usl{} $q$ belongs to $D(\vb)$. 
\rit{We finally claim that $\vb$ extends $q$}.   
This contradicts the contrary assumption above  
and completes the proof of \ref{l13}. 

To prove the extension claim, one has to check  
that $q(i)=b_i$ for all $i\in U$.
If $i\in U$ then $b_i=a_{\pi(i)}=a_{k_i}=q(i)$ 
by construction, as required. 

To prove claim \ref{l14} of the lemma, 
suppose otherwise.
This is forced by a \usl{} $p\in \qq W$, 
such that no function $\vb\in W^\om,$ 
\dd{\qq W}generic over $\rL(W)$ and 
extending $p$, is  
\dd{\dP\lom}generic over $\rL$. 
Arguing as in the proof of \ref{l13} above, 
we get a suitable permutation $\pi$ that yields  
a function $\vb\in W^\om,$ \dd{\qq W}generic 
over $\rL(W)$ and in the same time  
\dd{\dP\lom}generic over $\rL$ along with $\va$, 
and satisfies $W(\vb)=W(\va)=W$ 
(as a finite permutation of $\va$),  
and extends the \usl{} $p$. 
Therefore $\vb$ is \dd{\qq W}generic  
over $\rL(W)$ by claim \ref{l13} already established. 
This is a contradiction.

\ref{l15}
This is a generally known fact, yet we add a short proof.
As $Z\sq\rL(W)$,
there is a set $X\in\rL(W)$  
with $Z\sq X.$
Consider \dd{\qq W}names $s,t\in\rL(W)$ such that  
$Z=\bint{s}\va=\bint{t}\vb$, 
where $\bint{s}\va$ denotes the \dd\va interpretation 
of any given \dd{\qq W}name $s.$
By genericity, the equality $\bint{s}\va=\bint{t}\vb$
is forced by a pair of \usl s $p,q\in\qq W$, \ie\ 
$\va$ extends $p$, $\vb$ extends $q$, 
and if a pair $\ang{\va',\vb'}$ 
is \dd{(\qq W\ti\qq W)}generic over  
$\rL(W)$ and $\va'$ extends $p$, $\vb'$ extends $q$, 
then $\bint{s}{\va'}=\bint{t}{\vb'}$. 
\rit{We claim that the \usl{} $p$ \dd{\qq W}decides 
over $\rL(W)$ every  
sentence of the form $x\in\bint{s}{\dot\va}$}, where 
$\dot\va$ is a canonical \dd{\qq W}name for 
the principal generic function in $W^\om$. 

Indeed otherwise there exist functions 
$\va',\va''\in W^\om$, \dd{\qq W}generic over 
$\rL(W)$ and extending the \usl{} $p$, 
and an element $x\in X$, 
such that $x\in\bint{s}{\va'}$ but $x\nin\bint{s}{\va''}$. 
Consider a function $\vb'\in W^\om$, 
\dd{\qq W}generic both over $\rL(W)[\va']$ 
and over $\rL(W)[\va'']$, and extending the \usl{} $q$.
Then either pair $\ang{\va',\vb'}$, $\ang{\va'',\vb'}$ 
is \dd{(\qq W\ti\qq W)}generic over  
$\rL(W)$, but at least one of the two equalities  
$\bint{s}{\va'}=\bint{t}{\vb'}$,
$\bint{s}{\va''}=\bint{t}{\vb'}$
definitely fails, which is a contradiction. 

Thus  $p$ indeed \dd{\qq W}decides over $\rL(W)$ 
every sentence $x\in\bint{s}{\dot\va}$. 
This implies
$$
Z=\ens{x\in X}
{\text{$p$ \dd{\qq W}forces }x\in\bint{s}{\dot\va}
\text{ in }\rL(W)}
\in\rL(W)\,.
$$

\ref{l16}
To prove this key claim we apply a method  
introduced in \cite{kl31}. 
Consider a countable $\od$ set $Z\sq W^\om$ 
in $\rL[\va]$.
{\ubf Suppose towards the contrary that $Z\not\sq\rL(W)$}.

There is a formula $\vpi(z)$ 
with an ordinal $\ga_0$ as a parameter, 
such that $Z=\ens{z\in W^\om}{\vpi(z)}$ in $\rL[\va]$. 
There also exists a \usl{}
$p_0\in\qq W,$ $p_0\su\va$, 
which forces our assumptions, that is\vom
\ben 
\nenu
\itlb{enu1}%
$p_0$ \dd{\qq W}forces, over $\rL(W)$, 
that the set   
$\ens{z\in W^\om}{\vpi(z)}$ 
is countable and is not included in $\rL(W)$, 
or equivalently,   
if\/ $\vb\in W^\om$ is 
$\qq W$-generic over\/ $\rL(W)$ and  
extends $p_0$ 
then it is true in the extension 
$\rL(W)[\vb]=\rL[\vb]$
that the set  
$\Phi_\vb=\ens{z\in W^\om}{\vpi(z)}$ is countable and 
$\sus z\:(z\nin \rL(W)\land\vpi(z))$. 
It follows from the countability that there is a map 
$f_\vb:\om\na\Phi_\vb$, $f_\vb\in\rL[\vb]$.\vom
\een
Let $T\in\rL(W)$ be a canonical \dd{\qq W}name 
for $f_\vb$, so $f_\vb=\bint T\vb$.  
\vyk{
Further there exists a sequence 
$\sis{t_n}{n<\om}\in \rL(W)$ of \dd{\qq W}names, 
such that if a $\qq W$-gneric function  
$\vb\in W^\om$ extends $p_0$ then  
it is true in $\rL[\vb]=\rL(W)[\vb]$ that 
$$
\ens{z\in W^\om}{\vpi(z)}=\ens{\bint{t_n}\vb}{n<\om}\,, 
$$
where $\bint{t}\vb$ is the \dd\vb interpretation 
of the \dd{\qq W}name $t.$
Consider a canonical \dd{\qq W}name 
$T\in\rL(W)$ for the set   
$\ens{\bint{t_n}{\vb}}{n<\om}$.
}%
Then \ref{enu1} implies:\vom
\benq
\nenu
\atc
\itlb{enu2}%
$p_0$ \dd{\qq W}forces  
$\bint{\ran T}{\dva}=
\ens{\bint T{\dva}(n)}{n<\om}=  
\ens{z\in W^\om}{\vpi(z)} \not\sq\rL(W)$ 
over $\rL(W)$, 
or equivalently, if $\vb\in W^\om$   
$\qq W$-generic over $\rL(W)$ 
and $p_0\su \vb$ then it is true in $\rL[\vb]$
that  
$$\ran\bint{T}{\vb}=
\ens{\bint T{\vb}(n)}{n<\om}=  
\ens{z\in W^\om}{\vpi(z)}
\not\sq\rL(W).
$$ 
\eenq\vhm

Now our goal will be to  
{\ubf get a contradiction from \ref{enu2}}. 
Consider a regular uncountable cardinal $\ka>\ga_0$, 
such that the set $\rL_\ka$ is an elementary submodel 
of $\rL$ \poo\ a fragment of $\ZFC$ sufficiently 
large to prove the part of Lemma \ref{l1} already 
established including both \ref{enu1} and \ref{enu2}.
Then the set $\rL_\ka(W)$ contains $\ga_0$ and the   
name $T$.
As elements of the model $\rL_\ka(W)\sq\rL_\ka[\va]$, 
the sets $W,T$ admit canonical  
\dd{\dP\lom}names in $\rL_\ka$. 
Consider a countable elementary submodel 
$\gM\in\rL$ of $\rL_\ka$, 
containing those names and $\ga_0$. 
Then the sets $W,T$ 
and the forcing notion $\qq W$  
belong to $\gM(W)$. 
Consider the Mostowski collapse map 
$\pi:\gM(W)\na\rL_\la(W)$  
onto a transitive set of the form $\rL_\la(W)$, 
countable in $\rL[\va]$, where $\la<\omil$.
As $W$ is countable, we  have 
$\pi(W)=W$, $\pi(T)=T$, 
and hence $T\in\rL_\la(W)$, $\qq W\in\rL_\la(W)$.

We assert that there is $\vb\in W^\om$ satisfying\vhm 

\ben
\nenu
\atc
\atc
\itlb{enu3}\msur%
$\rL[\vb]=\rL[\va]$, 
$\vb$ is a \dd{\qq W}generic function over  
$\rL(W)$, $p_0\su\vb$, and the pair $\ang{\va,\vb}$ 
is \dd{(\qq W\ti\qq W)}generic over 
$\rL_\la(W)$.
\een

Indeed, as the set $\rL_\la$ is countable in  
$\rL$, there exists a bijection $h:\om\na\om,$ 
$h\in\rL$, equal to the identity on the (finite) 
domain $\dom{p_0}$ of the 
\usl{} $p_0\in\dC(W)$ (see above on $p_0$), 
and generic over $\rL_\la$  
in the sense of the Cohen-style  
forcing notion $\dB$ which consists of all 
injective tuples $u\in\om\lom.$ 
Let $\vb(n)=\va(h(n))$ for all $n$,  
\ie\ $\vb=\va\circ h$ is a superposition. 
Let's check that $\vb$ satisfies \ref{enu3}. 

Indeed, the function $\va$ of Lemma~\ref{l1} 
is generic over $\rL$, hence it is generic over  
$\rL_\la[h]\in\rL$, and hence the bijection $h$ 
is $\dB$-generic over $\rL_\la[\va]$ by the 
product forcing theorem.  
Therefore $h$ is generic over $\rL_\la(W)$, 
a smaller model. 
However $\va$ is \dd{\qq W}generic 
over $\rL_\la(W)$ by \ref{l13} of the lemma. 
It follows that the pair  
$\ang{\va,h}$ is 
\dd{(\qq W\ti\dB)}generic over $\rL_\la(W)$ 
still by the product forcing theorem. 
One easily proves then that 
$\ang{\va,\vb}$ is 
\dd{(\qq W\ti\qq W)}generic over $\rL_\la(W)$. 

We further have $\rL[\vb]=\rL[\va]$, because 
$h\in\rL$. 
Moreover  $\vb$ is \dd{\qq W}generic over 
$\rL(W)$, since $h\in\rL$ 
induces an order isomorphism   
of $\qq W$ in $\rL(W)$. 
Finally  $h$ is compatible with $p_0$ because $h$ 
is the identity on $\dom{p_0}$ 
by construction. 
This completes the proof that $\vb=\va\circ h$ 
satisfies \ref{enu3}.

In particular $\ww(\vb)=\ww(\va)=W$ holds, 
$\ran\bint{T}{\va}=\ran\bint{T}{\vb}\not\sq\rL(W)$ by 
\ref{enu2}.
On the other hand, the set  
$Z=\ran\bint{T}{\va}=\ran\bint{T}{\vb}$, 
belongs to the intersection  
$\rL_\la(W)[\va]\cap\rL_\la(W)[\vb]$ 
by construction. 
We conclude that $Z$ belongs to $\rL(W)$
by \ref{l15} of the lemma. 
(The above proof \ref{l15} is valid for 
$\rL_\la$ instead of $\rL$ as the ground model.)
The contradiction obtained completes the proof  
of \ref{l16}.\vom

\epF{Lemma \ref{l1} and Theorem \ref{t1}}

\punk{Comments and questions}
\las{x7}

Coming back to the Cohen-generic extensions, 
recall that if $a$ is a Cohen generic real over $\rL$ 
then $\HNT=\rL$ in $\rL[a]$ by Theorem~\ref{t3}. 

\bqe
\lam{QQ}
Is it true in generic extensions of $\rL$ by a single 
Cohen real that any countable \OD\ set consists of \OD\ 
elements?
\eqe

We cannot solve this even for \rit{finite} \OD\ sets. 
By the way it is not that obvious to expect the \rit{positive} 
answer. 
Indeed, the problem solves in the \rit{negative} for 
Sacks and some other generic extensions even for \rit{pairs}, 
see \cite{kl54,kl62}. 
For instance, if $a$ is a Sacks-generic real over $\rL$ then 
it is true in $\rL[a]$ that there is an $\OD$ unordered pair 
$\ans{X,Y}$ of sets of reals $X,Y\sq\pws\dn$ such that $X,Y$ 
themselves are non-\OD\ sets. 
See \cite{kl54} for a proof of this rather surprising 
result originally by Solovay.


\bibliographystyle{amsplain}

\bibliography{63,63klarx}

\end{document}